\documentclass[reqno]{amsart}
\usepackage{amscd}
\usepackage{amssymb}
\usepackage{amsmath}
\usepackage{amsthm}
\usepackage{latexsym}
\usepackage{upref}
\usepackage{epsfig}
\usepackage{mathrsfs}
\usepackage{float}
\usepackage[colorlinks,urlcolor=black]{hyperref}
\usepackage{graphics,graphicx,color}
\usepackage{floatflt}
\usepackage{indentfirst}
\usepackage[cmtip,arrow]{xy}
\usepackage{pb-diagram,pb-xy}
\usepackage{cite}

\newtheorem{thm}{Theorem}[section]
\newtheorem{prop}[thm]{Proposition}
\newtheorem{cor}[thm]{Corollary}
\newtheorem{lem}[thm]{Lemma}
\theoremstyle{definition}

\newtheorem{exam}[thm]{Example}
\numberwithin{equation}{section}

\newcommand{\vtu}{\ensuremath{\text{\textbf{\textit{u}}}}}
\newcommand{\vtv}{\ensuremath{\text{\textbf{\textit{v}}}}}

\newcommand{\vtx}{\ensuremath{\text{\textbf{\textit{x}}}}}
\newcommand{\vty}{\ensuremath{\text{\textbf{\textit{y}}}}}

\begin{document}
\title{on the normal centrosymmetric Nonnegative inverse eigenvalue problem}

\author{Somchai Somphotphisut}
\address{Department of Mathematics and Computer Science, Faculty of Science,
Chulalongkorn University, Phyathai Road, Patumwan, Bangkok 10330,
Thailand} \email{somchai.so@student.chula.ac.th \textrm{and
}kwiboonton@gmail.com}

\author{Keng Wiboonton}
\address{}
\email{}

\begin{abstract}
We give sufficient conditions of the nonnegative inverse eigenvalue problem (NIEP) for normal centrosymmetric matrices. These sufficient conditions are analogous to the sufficient conditions of the NIEP for normal matrices given by Xu \cite{xu} and Julio, Manzaneda and Soto \cite{julio2}.
\end{abstract}
\maketitle

\section{Introduction}
\medskip
The problem of finding necessary and sufficient conditions for $n$ complex numbers to be the spectrum of a nonnegative matrix is known as the \textit{nonnegative inverse eigenvalues problem} (NIEP). If the complex numbers $z_1, z_2, \ldots, z_n$ are the eigenvalues of an $n \times n$ nonnegative matrix $A$ then we say that the list $z_1, z_2, \ldots, z_n$ is realizable by $A$. The NIEP was first studied by Suleimanova in \cite{sul} and then several authors had been extensively studied this problem for examples in \cite{borobia, Kel, rado, salz, soto1, soto2, sul, loewy, reams, loewy2, wuwen}. 

There are related problems to the NIEP when one also want the desired nonnegative matrices to be symmetric, persymmetric, bisymmetric or circulant, etc. For example, if one want to find a necessary and sufficient condition for a list of $n$ numbers to be realized by a symmetric nonnegative matrix then this problem is called the NIEP for symmetric matrices or shortly the SNIEP (symmetric nonnegative inverse eigenvalue problem). In this paper, we are interested in the NIEP for normal centrosymmetric matrices.
 




There are several results of the NIEP for normal matrices. The following theorem due to Xu \cite{xu} gives a sufficient condition of the NIEP for normal matrices.
 

\bigskip
\begin{thm} \label{radwan1}
Let $\lambda_0 \geq 0 \geq \lambda_1 \geq \ldots \geq \lambda_n$ be real numbers and $z_1$, $z_2$, $\ldots$, $z_n$, $\overline{z}_1$, $\overline{z}_2$, $\ldots$, $\overline{z}_m$ be complex numbers such that \  $\textnormal{Im} \ z_j \neq 0$  for all $j$. If $\lambda_0 + \displaystyle \sum_{j=1}^n \lambda_j - 2\sum_{j=1}^m \vert z_j \vert \geq 0$ then there is an $(n+2m+1) \times (n+2m+1)$ nonnegative normal matrix with the eigenvalues $\lambda_0, \lambda_1,\ldots, \lambda_n, z_1, z_2, \ldots, z_n, \overline{z}_1, \overline{z}_2, \ldots, \overline{z}_m$.
\end{thm}
\bigskip

In \cite{julio2}, Julio, Manzaneda and Soto found another sufficient condition of the NIEP for normal matrices as the following.
\bigskip
\begin{thm} \textnormal{(Julio, Manzaneda, Soto  \cite{julio2}, 2015)} \label{julio5}
Let $\Lambda = \lbrace \lambda_1, \lambda_2, \ldots, \lambda_n \rbrace$ be a list of complex numbers with $\Lambda = \bar{\Lambda}$, $\lambda_1 \geq \textnormal{max}_j \vert \lambda_j \vert$ \ for $j= 2,3, \ldots, n,$ and $\displaystyle \sum_{j=1}^n \lambda_j \geq 0$ and let $\Lambda =\Lambda_0 \cup \Lambda_1 \cup \ldots \cup \Lambda_{{p_0}}$ be partition with

\qquad \qquad \qquad $\Lambda_0=\lbrace \lambda_{01}, \lambda_{02}, \ldots, \lambda_{0p_0} \rbrace$, \qquad $\lambda_{01} = \lambda_1$,

\qquad \qquad \qquad $\Lambda_k=\lbrace \lambda_{k1}, \lambda_{k2}, \ldots, \lambda_{kp_k} \rbrace$, \qquad $k=1, 2, \ldots, p_0$, \\
where some of the lists $\Lambda_k \ $ $(k=1,2,\ldots, p_0$ can be empty. Suppose that the following conditions hold:

\begin{itemize} 
\item[\textnormal{(1)}] For each $k=1, 2, \ldots, p_0$, there exists a normal nonnegative matrix with all eigenvalues as $\omega_k, \lambda_{k1}, \lambda_{k2}, \ldots, \lambda_{kp_k}$, $0\leq \omega_k \leq \lambda_1$.

\item[\textnormal{(2)}] There exists a normal nonnegative matrix $B$ of order $p_0$ with eigenvalues $\lambda_{01}, \lambda_{02}, \ldots, \lambda_{0p_0}$ and diagonal entries $\omega_1$ $\geq$ $\omega_2$ $\geq$ $\ldots$ $\geq$ $\omega_{p_0}$.
\end{itemize}
\noindent Then there is an $n \times n$ normal nonnegative matrix with the eigenvalues $\lambda_1, \lambda_2, \ldots, \lambda_n$.
\end{thm}
\bigskip
In this article we improve these results by adding some conditions so that the desired matrices are also centrosymmetric. Therefore we obtain sufficient conditions of the NIEP for normal centrosymmetric matrices.
The following is the summary of this article. In Section 2, we review basic properties of centrosymmetric matrices. Then we give a sufficient condition (analogous to the one in Theorem \ref{radwan1}) of the NIEP for normal centrosymmetric matrices in Section 3. The main result of Section 3 is Theorem \ref{centromain}. Finally, in Section 4, we provide a sufficient condition (analogous to the one in Theorem \ref{julio5}) of the NIEP for normal centrosymmetric matrices.   


\bigskip
\section{Centrosymmetric matrices}
\bigskip
Recall that an $n \times n$ matrix is called a centrosymmetric matrix if $JAJ=A$, where $J$ is the reverse identity matrix, i.e., $J=\begin{bmatrix}
 0 &&1\\
&\reflectbox{$\ddots$}&\\
1&& 0
\end{bmatrix}$.

 
The following well-known properties of centrosymmetric matrices which can be found in \cite{iyad} are collected below.

\begin{thm} \label{centroprop1}

Let $Q$ be an $n \times n$ centrosymmetric matrix.
\begin{itemize} 
\item[$\textnormal{(1)}$] If $n$ is an even number, then $Q$ is of the form
$\begin{bmatrix}
A & JCJ \\
C & JAJ
\end{bmatrix}$, where $A$ and $C$ are $\frac{n}{2} \times \frac{n}{2}$ matrices.

\item[$\textnormal{(2)}$] If $n=2m+1$ is an odd number, then $Q$ is of the form
$\begin{bmatrix}
A      & \vtx  & JCJ \\
\vty^T & p     &\vty^TJ \\
C      & J\vtx &JAJ
\end{bmatrix}$, where $\vtx$ and $\vty$ are $m \times 1$ matrices, $A$ and $C$ are $m \times m$ matrices.
\end{itemize}
\end{thm}

\bigskip

\begin{thm} \label{centroprop2}

Let $Q$ be an $n \times n$ centrosymmetric matrix.
\begin{itemize} 
\item[$\textnormal{(1)}$] If $n$ is an even number and $Q$ is of the form $Q=\begin{bmatrix}
A & JCJ \\
C & JAJ
\end{bmatrix}$, 
then $Q$ is orthogonally similar to the matrix
$\begin{bmatrix}
A-JC &  \\
 & A+JC
\end{bmatrix}$. Moreover, $Q$ is a normal matrix if and only if $A-JC$ and $A+JC$ are normal matrices and if $Q$ is a nonnegative with the Perron root $\lambda_0$ then $\lambda_0$ is the Perron root of $A+JC.$

\item[$\textnormal{(2)}$] If $n=2m+1$ is an odd number and $Q$ is of the form
$Q=\begin{bmatrix}
A   & \vtx  & JCJ \\
\vty^T & p  &\vty^TJ \\
C   & J\vtx &JAJ
\end{bmatrix}$,
then $Q$ is orthogonally similar to the matrices
$\begin{bmatrix}
A-JC &          &             \\
     &    p     & \sqrt{2}\vty^T \\
     &\sqrt{2}\vtx &  A+JC
\end{bmatrix}$ and 
$\begin{bmatrix}
A-JC &          &              \\
     &A+JC      & \sqrt{2}\vtx  \\
     &\sqrt{2}\vty^T &  p
\end{bmatrix}$. Moreover, $Q$ is a normal matrix if and only if $A-JC$ and 
$\begin{bmatrix}
p            & \sqrt{2}\vty^T \\
\sqrt{2}\vtx &  A+JC
\end{bmatrix}$ are normal matrices and if $Q$ is a nonnegative matrix with the Perron root $\lambda_0$, then $\lambda_0$ is the Perron root of 
$\begin{bmatrix}
p & \sqrt{2}\vty^{T} \\ 
\sqrt{2}\vtx & A+JC 
\end{bmatrix}$.
\end{itemize}
\end{thm}

\medskip

\medskip
Next, we give a trivial fact concerning the eigenvector corresponding to the Perron root of a nonnegative centrosymmetric matrix.
\medskip
\begin{prop} \label{centroprop4}
Let $Q$ be an $n \times n$ nonnegative centrosymmetric matrix. If $\vtv$ is an eigenvector corresponding to an eigenvalue $\lambda$ of $Q$ then $J\vtv$ is also an eigenvector corresponding to the eigenvalue $\lambda$. Moreover, if $\lambda_0$ is the Perron root of $Q$ then there is a nonnegative eigenvector $\vtv_0$ such that $J\vtv_0=\vtv_0$.
\end{prop}

\bigskip
\section{A Sufficient condition of the NIEP for centrosymmetric matrices}
\medskip
In this section, we will investigate the NIEP for normal centrosymmetric matrices. So, we have the spectral theorem for normal matrices as a tool for our investigation. 


Note that since any normal centrosymmetric matrix with all real entries and real eigenvalues is symmetric and any symmetric centrosymmetric matrix is bisymmetric, for a list of $n$ real numbers the problem of the NIEP for normal centrosymmetric matrices can be solved by using the results of the NIEP for bisymmetric matrices which can be found in \cite{julio}. So, from now on we study the NIEP for normal centrosymmetric matrices with a list of $n$ complex numbers not all of them are real.


For $n \leq 3$, we have the following result of the NIEP for centrosymmetric matrices.

\begin{thm}
Let $Q$ be an $n \times n$ nonnegative centrosymmetric matrix with $n \leq 3$. Then all eigenvalues of $Q$ are real. Moreover, $\lambda_1 \geq \ldots \geq \lambda_n$ are the eigenvalues of $Q$ if and only if $\lambda_1 \geq \vert \lambda_n \vert$ and $\sum_{j=1}^{n} \lambda_{j} \geq 0$.
\begin{proof}
The assertion is trivial for $n=1$. If $n=2$ then any $2 \times 2$ centrosymmetric matrix is bisymmetric and hence its eigenvalues are real. Moreover, it is easy to show that $\lambda_1 \geq \lambda_2$ are eigenvalues of a $2 \times 2$ bisymmetirc nonnegative matrix if and only if $\lambda_1 \geq \vert \lambda_2 \vert$. 

Now, we prove the assertion for $n=3$.
Write $Q=\begin{bmatrix}
\alpha & \beta & \gamma \\
\xi    & p     &  \xi   \\
\gamma & \beta & \alpha
\end{bmatrix}$. 
By Theorem \ref{centroprop2}, $Q$ is orthogonally similar to 
$\begin{bmatrix}
\alpha-\gamma  &                 &                \\
               &   p             & \sqrt{2}\xi \\
               &  \sqrt{2}\beta  & \alpha + \gamma
\end{bmatrix}$.
Since $\begin{bmatrix}
p            & \sqrt{2}\xi \\
\sqrt{2}\beta  & \alpha + \gamma
\end{bmatrix}$
is a nonnegative matrix, its eigenvalues must be real numbers by the Perron-Frobenius theorem and the fact that the complex eigenvalues of a nonnegative matrix must be in conjugate pairs. This implies that the eigenvalues of $Q$ are real. If $\lambda_1 \geq \lambda_2 \geq \lambda_3$ are all eigenvalues of $Q$ then $\lambda_1 \geq \vert \lambda_3 \vert $ and $\sum_{j=1}^3 \lambda_j \geq 0$ by the well-known necessary condition of the NIEP. 

Conversely, let $\lambda_1 \geq \lambda_2 \geq \lambda_3$ be real numbers. The condition $\lambda_1 \geq \vert \lambda_3 \vert $ and $\sum_{j=1}^3 \lambda_j \geq 0$ is a sufficient condition for an existence of a $3 \times 3$ bisymmetric nonnegative matrix with the eigenvalues $\lambda_1, \lambda_2, \lambda_3$ (see \cite{julio}) and hence a sufficient condition for the case of centrosymmetric matrices.  
\end{proof}
\end{thm}
\medskip
Next, we consider the NIEP for normal centrosymmetric matrices when $n=4$. In this case, we have a necessary and sufficient condition as the following theorem.
\medskip
\begin{thm} \label{centro1}
Let $\lambda_0$ and $\lambda_1$ be real numbers and $z$ = $a+bi$ where $a \in \mathbb{R}$ and $b > 0$. Suppose that $\lambda_0 \geq \textnormal{max} \ \lbrace \lambda_1, \vert z \vert \rbrace$. Then $\lambda_0, \lambda_1, z, \overline{z}$ are the eigenvalues of a $4 \times 4$ normal centrosymmetric nonnegative matrix if and only if $\lambda_0+\lambda_1 - 2 \vert a \vert \geq 0$ and $\lambda_0 - \lambda_1 - 2\vert b \vert \geq 0$. 

\begin{proof}
Suppose that $\lambda_0, \lambda_1, z, \overline{z}$ are the eigenvalues of a $4 \times 4$ normal centrosymmetric nonnegative matrix 
$Q=\begin{bmatrix}
A&JCJ\\
C&JAJ
\end{bmatrix},$ where $A$ and $C$ are $2 \times 2$ nonnegative matrices.
Then, by Theorem \ref{centroprop2}, $\lambda_0$ and $\lambda_1$ must be the eigenvalues of the $2 \times 2$ normal matrix $A+JC$ and $z, \overline{z}$ are the eigenvalues of the $2 \times 2$ normal matrix $A-JC$. 

Since $A+JC$ is a $2 \times 2$ nonnegative normal matrix with real eigenvalues, $A+JC$ must be a nonnegative symmetric  matrix.
Write $A+JC$ as 
$\begin{bmatrix}
\alpha&\beta\\
\beta&\gamma
\end{bmatrix},$ where $\alpha, \beta, \gamma \in \mathbb{R}.$
Note that any $2 \times 2$ real normal matrix with the complex eigenvalues $z$ and $\overline{z}$ is either the matrix 
$\begin{bmatrix}
a&-b\\
b&a
\end{bmatrix}$ or the matrix
$\begin{bmatrix}
a &b\\
-b&a
\end{bmatrix}$.
Without loss of generality, we may assume that $A-JC=\begin{bmatrix}
a&-b\\
b&a
\end{bmatrix}.$
Since $2A=(A+JC) + (A-JC)$ and $JC=(A+JC)-(A-JC)$ are nonnegative matrices, we have $\alpha \geq \vert a \vert, \beta \geq \vert b \vert$ and $\gamma \geq \vert a \vert.$
By comparing the characteristic polynomial of $A+JC$ with $(x-\lambda_0)(x-\lambda_1)$, we have $\alpha + \gamma = \lambda_0+\lambda_1$ and $\alpha\gamma - \beta^2 = \lambda_0\lambda_1.$ 
Then $\lambda_0+\lambda_1 = \alpha + \gamma \geq \vert a \vert  + \vert a \vert = 2 \vert a \vert.$ Thus $\lambda_0 + \lambda_1 - 2\vert a \vert \geq 0$.

Next, note that $\lambda_0 \lambda_1 + \beta^2 = \alpha \gamma = \alpha (\lambda_0 + \lambda_1 - \alpha).$
This implies that $\alpha^2 - (\lambda_0 + \lambda_1)\alpha + \lambda_0\lambda_1 + \beta^2=0$ and hence $\alpha$ is a real solution of the quadratic equation $x^2 - (\lambda_0 + \lambda_1)x + \lambda_0\lambda_1 + \beta^2=0$. Thus $(\lambda_0+\lambda_1)^2 - 4(\lambda_0\lambda_1+\beta^2) \geq 0$.
So, $(\lambda_0 - \lambda_1-2\beta)(\lambda_0 - \lambda_1+2\beta) = (\lambda_0 - \lambda_1)^2 - 4\beta^2 = (\lambda_0+\lambda_1)^2 - 4(\lambda_0\lambda_1+\beta^2) \geq 0$.
Since $(\lambda_0 - \lambda_1+2\beta) \geq 0$, $(\lambda_0 - \lambda_1-2\beta)\geq 0.$ This implies that $\lambda_0 - \lambda_1-2 \vert b \vert \geq 0$ because $\beta \geq  \vert b \vert$.

Conversely, if we have $\lambda_0 + \lambda_1-2 \vert a \vert \geq 0$ and $\lambda_0 - \lambda_1-2 \vert b \vert \geq 0,$ then the matrix
$Q=\dfrac{1}{4}\begin{bmatrix}
A&JCJ\\
C&JAJ
\end{bmatrix}$, where $A=\begin{bmatrix}
\lambda_0+\lambda_1+2a&\lambda_0-\lambda_1-2b\\
\lambda_0-\lambda_1+2b&\lambda_0+\lambda_1+2a
\end{bmatrix}$ and 
$C=\begin{bmatrix}
\lambda_0-\lambda_1-2b&\lambda_0+\lambda_1-2a\\
\lambda_0+\lambda_1-2a&\lambda_0-\lambda_1+2b
\end{bmatrix}$, is 
a $4 \times 4$ normal centrosymmetric nonnegative matrix with the eigenvalues $\lambda_0, \lambda_1, z, \overline{z}.$
\end{proof}
\end{thm}
\bigskip
\begin{exam}
The list of four complex numbers $10, 5, 3 \pm 4i$ satisfies the condition of Theorem \ref{julio5}. Moreover, this list satisfies the condition of Theorem $3.3$ in \cite{nizar}. Then there is a $4 \times 4$ normal nonnegative matrix with the eigenvalues $10, 5, 3 \pm 4i$. However, the list $10, 5, 3+4i$ does not satisfy the condition of Theorem \ref{centro1}, therefore there is no $4 \times 4$ normal centrosymmetric nonnegative matrix with the eigenvalues $10, 5, 3 \pm 4i$.
\end{exam}
\bigskip
The following theorem due to Xu \cite{xu} is a result for combining two normal nonnegative matrices with a perturbation on the Perron roots . 
\bigskip
\begin{thm} \label{normal1}
Let $A$ be an $(m+1) \times (m+1)$ nonnegative normal matrix with the eigenvalues $\alpha_0, \alpha_1, \ldots ,\alpha_m$ and $\alpha_0 \geq \vert \alpha_j \vert$, $j=1,2, \ldots,m$ and $\vtu$ be a unit nonnegative eigenvector corresponding to $\alpha_0 \textnormal{;}$ let $B$ be an $(n+1) \times (n+1)$ nonnegative normal matrix with the eigenvalues $\beta_0$, $\beta_1$, $\ldots$ , $\beta_n$ and $\beta_0 \geq \vert \beta_j \vert$, $j=1,2, \ldots, n$ and $\vtv$ be a unit nonnegative eigenvector corresponding to $\beta_0$ and $\alpha_0 \geq \beta_0.$ Then for any nonnegative number $\rho$, the matrix 
$$\begin{bmatrix}
A & \rho \vtu \vtv^T \\
\rho \vtv \vtu^T & B
\end{bmatrix}$$ is an $(m+n+2) \times (m+n+2)$ nonnegative normal matrix with the eigenvalues $\gamma_1$, $\gamma_2$, $\alpha_1$, $\alpha_2$, $\ldots$, $\alpha_m$, $\beta_1$, $\beta_2$, $\ldots$, $\beta_n$, where $\gamma_1, \gamma_2$ are the eigenvalues of the matrix $$\begin{bmatrix}
\alpha_0 & \rho \\
\rho     &  \beta_0
\end{bmatrix}.$$
\end{thm}
\bigskip
For our purpose, we mimic the above theorem of Xu to the case of normal centrosymmetric matrices. Our result is the following lemma.
\bigskip
\begin{lem} \label{ourcc2} 
Let $A$ be an $m \times m$ normal centrosymmetric nonnegative matrix with the eigenvalues $\alpha_1, \alpha_2, \ldots, \alpha_m$ such that $\alpha_1$ is the Perron root of $A$ and let $\vtu_1$ be a unit nonnegative eigenvector corresponding to $\alpha_1$ such that $J\vtu_1=\vtu_1 \textnormal{;}$ let $B$ be an $n \times n$ nonnegative normal centrosymmetric matrix with the eigenvalues $\beta_1, \beta_2, \ldots, \beta_n$ such that $\beta_1$ is the Perron root of $B$ and let $\vtv_1$ be a unit nonnegative eigenvector corresponding to $ \beta_1$ such that $J\vtv_1=\vtv_1$. If $\gamma_1 , \gamma_2 , \gamma_3$ are all eigenvalues of the matrix
$$\widehat{C} = 
\begin{bmatrix}
\beta_1 & \rho & \xi \\ \rho & \alpha_1 & \rho \\ \xi & \rho & \beta_1
\end{bmatrix},$$
where $\rho$, $\xi \geq$ 0, then the matrix
$$C = \begin{bmatrix}
B                           &\rho \vtv_1 \vtu_1^{T}     &\xi \vtv_1 \vtv_1^{T} \\ 
\rho \vtu_1 \vtv_1^{T}      &      A                    & J\rho \vtu_1 \vtv_1^{T}J \\ 
\xi \vtv_1 \vtv_1^{T}       &\rho \vtv_1 \vtu_1^{T}     &B
\end{bmatrix}$$
is a normal centrosymmetric nonnegative matrix with all eigenvalues as
$$\gamma_1, \gamma_2, \gamma_3, \alpha_2, \alpha_3, \ldots, \alpha_m, \beta_2, \beta_3, \ldots, \beta_n, \beta_2, \beta_3 , \ldots, \beta_n.$$
\begin{proof}
Let $\vtu_1,\vtu_2, \ldots, \vtu_m$ and $\vtv_1$, $\vtv_2$, $\ldots$, $\vtv_n$ be orthonormal systems of the eigenvectors of $A$ and $B$ respectively and let 
$\begin{bmatrix} r_1 \\ s_1 \\ t_1 \end{bmatrix}, 
\begin{bmatrix} r_2 \\ s_2 \\ t_2 \end{bmatrix},
\begin{bmatrix} r_3 \\ s_3 \\ t_3 \end{bmatrix}$ 
be an orthonormal system of the eigenvectors of $\widehat{C}.$ Then 
$\begin{bmatrix} \textbf{0} \\ \vtu_j \\ \textbf{0} \end{bmatrix}, 
\begin{bmatrix} \vtv_k \\ \textbf{0} \\ \textbf{0} \end{bmatrix}, 
\begin{bmatrix} \textbf{0} \\ \textbf{0} \\ \vtv_k \end{bmatrix}, 
\begin{bmatrix} r_l \vtv_1 \\ s_l \vtu_1 \\ t_l \vtv_1 \end{bmatrix}$ 
are $m+2n$ orthogonal eigenvectors of $C$ corresponding to the eigenvalues $\alpha_j , \beta_k , \gamma_l$ for $j=2,\ldots, m, \ k=2,\ldots, n,$ and $l=1,2,3$, and hence by normalizing these vectors we obtain an orthonormal system for $C$. By the spectral theorem for normal matrices, $C$ is a normal matrix. We see that $C$ is nonnegative since $A$, $B$, $\vtu_1$, $\vtv_1$, $\rho$ and $\xi$ are nonnegative. Moreover, $C$ is centrosymmetric because $A$ and $B$ are centrosymmetric, $J\vtu_1=\vtu_1$ and $J\vtv_1=\vtv_1.$ 
\end{proof}
\end{lem}
\medskip

Now, if we apply Lemma \ref{ourcc2} (assume that $\alpha_1 \geq \beta_1$) by choosing $$\rho=\sqrt{\frac{-(\alpha_1-\beta_1-a)(a+b)}{2}} \ \text{ \ and \ } \xi = -b,$$
where $a$ and $b$ are real numbers such that $\alpha_1-\beta_1 \geq a \geq b$ and $a+b \leq 0$, then  we have that $\alpha_1 - (a+b)$, $\beta_1+a$ and $\beta_1+b$ are all eigenvalues of $\widehat{C}$. Therefore we have the following result.

\medskip
\begin{cor} \label{ourcc3}
Let $A$ and $B$ be normal centrosymmetric nonnegative matrices as in Lemma \ref{ourcc2} with the Perron root $\alpha_1$ and $\beta_1$, respectively with $\alpha_1 \geq \beta_1$. If $a$, $b$ are real numbers such that $\alpha_1-\beta_1 \geq a \geq b$ and $a+b \leq 0$ then there is a normal centrosymmetric nonnegative matrix with all eigenvalues as $\alpha_1 - (a+b), \beta_1+a$, $\beta_1+b$, $\alpha_2$ , $\alpha_3$, $\ldots$, $\alpha_m$ , $\beta_2$, $\beta_3$, $\ldots$, $\beta_n$, $\beta_2$, $\beta_3$, $\ldots$, $\beta_n$.
\end{cor}
\medskip

We now give a sufficient condition of the NIEP for centrosymmetric normal nonnegative matrices when there are at least two real numbers and only one pair of complex numbers in the list of given numbers.

\medskip
\begin{thm} \label{centro2}
Let $\lambda_0 \geq 0 > \lambda_1 \geq \ldots \geq \lambda_n$ be real numbers and $z = a+bi$ where $a \in \mathbb{R}$ and $b > 0$. If $\lambda_0 + \sum_{j=1}^n \lambda_j - 2\vert z \vert \geq 0$, then there is an $(n+3) \times (n+3)$ normal centrosymmetric nonnegative matrix with the eigenvalues $\lambda_0, \lambda_1, \ldots, \lambda_n, z, \overline{z}$.
\begin{proof}
We will prove this theorem by induction on $n$.
If $n=1$, then the assertion follows from Theorem \ref{centro1}.

Suppose that $n=2$. Then let $A$ and $C$ be nonnegative matrices such that $A-JC =
\begin{bmatrix}
a&-b\\
b&a
\end{bmatrix}$ and
$A+JC=\dfrac{1}{2}\begin{bmatrix}
\lambda_0 +\lambda_1+\lambda_2 &\lambda_0+\lambda_1-\lambda_2\\
\lambda_0 +\lambda_1-\lambda_2 &\lambda_0+\lambda_1+\lambda_2\\
\end{bmatrix}$
and let
$B=\begin{bmatrix}
A+JC&\sqrt{-\lambda_0 \lambda_1}\vtu\\
\sqrt{-\lambda_0 \lambda_1}\vtu^T&0
\end{bmatrix}$ where
$\vtu=\begin{bmatrix}
\dfrac{1}{\sqrt{2}}\\ 
\dfrac{1}{\sqrt{2}}\\
\end{bmatrix}$.
Then $A-JC$ and $B$ are normal matrices with the eigenvalues $z, \overline{z}$ and $\lambda_0, \lambda_1, \lambda_2$, respectively. So, the normal centrosymmetric nonnegative matrix 
$Q=\begin{bmatrix}
A   					  & \dfrac{\rho \vtu}{\sqrt{2}}   & JCJ                   \\
\dfrac{\rho \vtu^T}{\sqrt{2}}     & 0                          & \dfrac{\rho \vtu^TJ}{\sqrt{2}} \\
C                         & \dfrac{\rho J\vtu}{\sqrt{2}}  & JAJ
\end{bmatrix}$, where $\rho =\sqrt{-\lambda_0 \lambda_1}$,
has $\lambda_0, \lambda_1, \lambda_2, z, \overline{z}$ as its eigenvalues.

Next, let $n \geq 3$ and suppose that the assertion is true for all smaller systems of numbers satisfying the condition.
Let $\lambda^{\prime}_0=\lambda_0+\lambda_1+\lambda_2$. Then the system $\lambda^{\prime}_0, \lambda_3, \ldots, \lambda_n, z$ satisfies the condition and so by the induction hypothesis, there is an $(n+1) \times (n+1)$ normal centrosymmetric nonnegative matrix $\tilde{Q}$ with the eigenvalues $\lambda^{\prime}_0, \lambda_3, \ldots, \lambda_n, z, \overline{z}$. Applying Corollary \ref{ourcc3} with the matrices $\tilde{Q}$ and $\begin{bmatrix}
0
\end{bmatrix}$ and $a=\lambda_1, b=\lambda_2$, we have that there is an $(n+3) \times (n+3)$ normal centrosymmetric nonnegative matrix with the eigenvalues $\lambda_0, \lambda_1, \ldots, \lambda_n, z, \overline{z}$.   
\end{proof}
\end{thm}
\medskip

Next, we will give a sufficient condition of the NIEP for normal centrosymmetric matrices when there are only two real numbers in our list of complex numbers. Before we obtain this result we need the following two lemmas.

\medskip
\begin{lem} \label{tool1}
Let $\lambda$ be a real number and $z=a+bi$ where $a, b \in \mathbb{R}$. If $\lambda_0 \geq 2\vert z \vert \geq 0$, then $\lambda_0 \geq \vert a \vert +\sqrt{3} \vert b \vert.$
\begin{proof}
Note that $9a^4+6a^2b^2+b^4 \geq 12a^2b^2$ since $9a^4-6a^2b^2+b^4 = (3a^2-b^2)^2 \geq 0.$ It implies that $(3a^2+b^2)^2 \geq 12a^2b^2$ and hence $3a^2+b^2 \geq 2\sqrt{3} \vert{ab} \vert$. So, $4a^2+4b^2 \geq a^2+2\sqrt{3}\vert ab \vert + 3b^2.$ Thus, $2\sqrt{a^2+b^2} \geq \vert a \vert + \sqrt{3} \vert b \vert$. Therefore, if $\lambda_0 \geq 2 \vert z \vert$, then $\lambda_0 \geq 2\vert z \vert = 2\sqrt{a^2+b^2} \geq \vert a \vert +\sqrt{3}\vert b \vert.$
\end{proof}
\end{lem}
\bigskip
\begin{lem} \label{cc3}
Let $\lambda_0$ be a real number and $z_1, z_2, \ldots, z_m$ be complex numbers such that $z_j=a_j+b_ji$ where $a_j, b_j \in \mathbb{R}$ for $j=1, \ldots, m$. If $\lambda_0 - 2\sum_{j=1}^m \vert z_j \vert \geq 0$, then 
$$c_k :=\lambda_0 + 2\sum_{j=1}^m a_j \textnormal{cos} \left( \dfrac{2kj\pi}{2m+1} \right) +\displaystyle 2\sum_{j=1}^m b_j \textnormal{sin} \left( \dfrac{2kj\pi}{2m+1} \right) \geq 0$$
for $k=0, 1, \ldots, 2m.$
\begin{proof}
We first recall the following well-known inequality: 
$$\vert a \text{cos} \ \theta + b \text{sin} \ \theta \vert \ \leq \sqrt{a^2+b^2} \ \ \text{for} \ a, b, \theta \in \mathbb{R}.$$
Now, suppose that $\lambda_0 - 2\sum_{j=1}^m \vert z_j \vert \geq 0$. Then the above inequality implies that 
$$c_k:=\lambda_0 + 2\displaystyle \sum_{j=1}^m a_j \text{cos} \left( \dfrac{2kj\pi}{2m+1} \right) +\displaystyle 2\sum_{j=1}^m b_j \text{sin} \left( \dfrac{2kj\pi}{2m+1} \right) \geq 0$$
for $k=0, 1, \ldots, 2m.$
\end{proof}
\end{lem}
\medskip
Now, we ready to give a sufficient condition of the NIEP for normal centrosymmetric matrices when there are $4s$ complex numbers and only two real numbers in the given list of numbers. We employ the method in \cite{rojo2} which gives circulant matrices having the eigenvalues as the given list of complex numbers.
\medskip
\begin{thm} \label{centro3} 
Let $\lambda_0 \geq 0 > \lambda_1$ be real numbers and 
$z_1, z_2, \ldots, z_{2s}$ be complex numbers such that $z_j=a_j+b_ji$ where $a_{j} \in \mathbb{R}$ and $b_j >0$ for $j=1, 2, \ldots, 2s$.
If $\lambda_0 + \lambda_1 - 2 \sum_{j=1}^{2s} \vert z_j \vert \geq 0$, then there is a 
$(4s+2) \times (4s+2)$ normal centrosymmetric nonnegative matrix with the eigenvalues $\lambda_0, \lambda_1, z_1, \ldots, z_{2s}, \overline{z}_1, \ldots, \overline{z}_{2s}$.
\begin{proof}
Suppose that $\lambda_0 + \lambda_1 - 2\displaystyle \sum_{j=1}^{2s} \vert z_j \vert \geq 0$. Let $A$ and $C$ be $(2s+1) \times (2s+1)$ matrices such that 
$$A+JC= \begin{bmatrix}
c_0   & c_1    & c_2    & \cdots & c_{2s}  \\
c_{2s} & c_0    & c_1    & \cdots & c_{2s-1}\\
\vdots& \vdots & \vdots & \ddots & \vdots \\
c_1   & c_2    & c_3    & \cdots & c_0
\end{bmatrix} \text{and } A-JC=
\begin{bmatrix}
d_0   & d_1    & d_2    & \cdots & d_{2s}  \\
d_{2s} & d_0    & d_1    & \cdots & d_{2s-1}\\
\vdots& \vdots & \vdots & \ddots & \vdots \\
d_1   & d_2    & d_3    & \cdots & d_0
\end{bmatrix},$$ where
$$c_k=\dfrac{1}{2s+1} \begin{pmatrix}\lambda_0 + 2\displaystyle \sum_{j=1}^s a_j \text{cos} \left( \dfrac{2kj\pi}{2s+1} \right) +\displaystyle 2\sum_{j=1}^s b_j \text{sin} \left( \dfrac{2kj\pi}{2s+1} \right) \end{pmatrix}$$ 
and
$$d_k=\dfrac{1}{2s+1} \begin{pmatrix}\lambda_1 + 2\displaystyle \sum_{j=1}^s a_{s+j} \text{cos}\left(  \dfrac{2kj\pi}{2s+1} \right) + \displaystyle 2\sum_{j=1}^s b_{s+j} \text{sin} \left( \dfrac{2kj\pi}{2s+1} \right) \end{pmatrix}.$$ 
By the proof of Theorem $10$ in \cite{rojo2}, the eigenvalues of $A+JC$ and $A-JC$ 
are $\lambda_0, z_1, \ldots, z_s, \overline{z}_1, \ldots, \overline{z}_s$ and $\lambda_1, z_{s+1}, \ldots, z_{2s}$, $\overline{z}_{s+1}, \ldots, \overline{z}_{2s}$, respectively. 
Note that since $\lambda_0 + \lambda_1 - 2\displaystyle \sum_{j=1}^{2s} \vert z_j \vert \geq 0$, by Lemma \ref{cc3} and $\lambda_1 < 0$, $A+JC$ is nonnegative. Moreover, each entry of $A+JC$ is greater than or equal to the absolute value of the corresponding entry of $A-JC$ by Lemma \ref{cc3}. So, the matrix $Q = \begin{bmatrix}
A & JCJ \\
C & JAJ
\end{bmatrix}$
is a nonnegative centrosymmetric matrix with the eigenvalues $\lambda_0, \lambda_1, z_1, \ldots, z_{2s}, \overline{z}_1, \ldots, \overline{z}_{2s}$. Since any circulant matrix is a normal matrix, $A+JC$ and $A-JC$ are normal matrices. Therefore $Q$ is also a normal matrix and hence we are done.
\end{proof}
\end{thm}
\medskip
Next, we extend Theorem \ref{centro3} to the case that there are $n+1$ real numbers in the list of $(n+4s+1)$ numbers.
\medskip
\begin{thm} \label{centro4}
Let $\lambda_0 \geq 0 > \lambda_1 \geq \ldots \geq \lambda_n, n \geq 1,$ be real numbers and 
$z_1, z_2, \ldots, z_{2s}$ be complex numbers such that $z_j=a_j+b_ji$ where $a_{j} \in \mathbb{R}$ and $b_j >0$ for $j=1, 2, \ldots, 2s$.
If $\lambda_0 + \sum_{j=1}^n \lambda_j - 2 \sum_{j=1}^{2s} \vert z_j \vert \geq 0$, then there is an 
$(n+4s+1) \times (n+4s+1)$ normal centrosymmetric nonnegative matrix with all eigenvalues as $\lambda_0, \lambda_1, \ldots, \lambda_n, z_1, \ldots, z_{2s}, \overline{z}_1, \ldots, \overline{z}_{2s}$.
\begin{proof}
Suppose that $\lambda_0 + \sum_{j=1}^n \lambda_j -2 \sum_{j=1}^{2s} \vert z_j \vert \geq 0$. If $n=1$ then the assertion follows from Theorem \ref{centro3}.

Now, suppose that $n=2$ and let $\lambda^{\prime}_0=\lambda_0+\lambda_2$. So, the list $\lambda^{\prime}_0, \lambda_1, z_1, \ldots, z_s$ satisfies the inequality condition in Theorem \ref{centro3}. Following the proof of Theorem \ref{centro3}, we let 

$$A+JC= \begin{bmatrix}
c_0   & c_2    & c_3    & \cdots & c_{2s}  \\
c_{2s} & c_0    & c_2    & \cdots & c_{2s-1}\\
\vdots& \vdots & \vdots & \ddots & \vdots \\
c_1   & c_2    & c_3    & \cdots & c_0
\end{bmatrix},$$ 
where
$c_k=\dfrac{1}{2s+1} \begin{pmatrix}\lambda^{\prime}_0 + 2\displaystyle \sum_{j=1}^s a_j$ cos $\left( \dfrac{2kj\pi}{2s+1} \right)+\displaystyle 2\sum_{j=1}^s b_j$ sin $\left( \dfrac{2kj\pi}{2s+1} \right) \end{pmatrix}$.
We now extend the matrix $A+JC$ to the matrix 
$E=\begin{bmatrix}
A+JC                        & \sqrt{-\lambda_0 \lambda_2} \vtu \\
\sqrt{-\lambda_0 \lambda_2}\vtu^T & 0
\end{bmatrix}$, where $\vtu$ is a nonnegative unit eigenvector of $A+JC$ corresponding to $\lambda^{\prime}_0$. Then the matrix $E$ is a nonnegative normal matrix with the eigenvalues $\lambda_0, \lambda_2, z_1, \ldots, z_s$, $\overline{z}_1, \ldots, \overline{z}_{s}$ by Theorem \ref{normal1}. By defining $A-JC$ as in the proof of Theorem \ref{centro3}, we have that the matrix 
$Q=\begin{bmatrix}
A   					  & \dfrac{\rho \vtu}{\sqrt{2}}   & JCJ                   \\
\dfrac{\rho \vtu^T}{\sqrt{2}}     & 0                     & \dfrac{\rho \vtu^TJ}{\sqrt{2}} \\
C                         & \dfrac{\rho J\vtu}{\sqrt{2}}  & JAJ
\end{bmatrix}$, where $\rho = \sqrt{-\lambda_0 \lambda_2}$, is our desired matrix.

Now, we proceed by induction on $n \geq 3$. So assume that the assertion is true for all systems of $\lambda$'s having length less than $n$ and satisfying the condition in the assertion. Define $\lambda^{\prime}_0 = \lambda_0+\lambda_{n}+\lambda_{n-1}$. Then the list of numbers $\lambda^{\prime}_0, \lambda_1, \ldots, \lambda_{n-2}$, $z_1, \ldots, z_{2s}$ satisfies the condition in the induction hypothesis. Therefore there is a $(4s+n-1) \times (4s+n-1)$ normal centrosymmetric nonnegative matrix $Q$ with the eigenvalues $\lambda^{\prime}_0, \lambda_1, \ldots, \lambda_{n-2}$, $z_1, \ldots, z_{2s}, \overline{z}_1, \ldots, \overline{z}_{2s}.$ Applying Corollary \ref{ourcc3} to the matrices $Q$ and $\begin{bmatrix}
0
\end{bmatrix}$ with the real numbers $\lambda_{n-1}$ and $\lambda_{n}$, we have that there is a $(4s+n+1) \times (4s+n+1)$ normal centrosymmetric nonnegative matrix with the eigenvalues $\lambda_0, \lambda_1, \ldots, \lambda_n, z_1, \ldots, z_{2s}, \overline{z}_1, \ldots, \overline{z}_{2s}$.

\end{proof}
\end{thm}
\medskip
Next, we consider the case $n=8$ with only $2$ real numbers in our list.
\medskip
\begin{thm} \label{centro9}
Let $\lambda_0 \geq 0 > \lambda_1$ be real numbers and $z_1, z_2, z_3$ be complex numbers such that $z_j=a_j+b_{j}i$ where $a_{j} \in \mathbb{R}$ and $b_{j}>0$ for $j=1,2,3$. If $\lambda_0 + \lambda_1 - 2\sum_{j=1}^3 \vert z_{j} \vert \geq 0$ then there is an $8 \times 8$ normal centrosymmetric nonnegative matrix with the eigenvalues $\lambda_0, \lambda_1, z_1, z_2, z_3, \overline{z}_1, \overline{z}_2, \overline{z}_3.$
\begin{proof}
We first find matrices $A$ and $C$ such that

$$A-JC= 
\dfrac{1}{2}\begin{bmatrix}
a_2+a_3  &  -(b_2+b_3)   &  -(b_2-b_3)   & a_2 - a_3 \\
b_2+b_3  &   a_2 + a_3   &  a_2 - a_3    & b_2 - b_3  \\
b_2-b_3  &   a_2 -a_3    &  a_2+a_3      & b_2 + b_3  \\
a_2-a_3  &   -(b_2-b_3)  &  -(b_2+b_3)   & a_2 + a_3 
\end{bmatrix}$$ 
\noindent and 
$$A+JC=
\dfrac{1}{4}\begin{bmatrix}
\lambda_0+\lambda_1+2a_1  & \lambda_0-\lambda_1-2b_1  & \lambda_0-\lambda_1+2b_1  & \lambda_0+\lambda_1-2a_1    \\
\lambda_0-\lambda_1+2b_1  & \lambda_0+\lambda_1+2a_1  & \lambda_0+\lambda_1-2a_1  & \lambda_0-\lambda_1-2b_1    \\
\lambda_0-\lambda_1-2b_1  & \lambda_0+\lambda_1-2a_1  & \lambda_0+\lambda_1+2a_1  & \lambda_0-\lambda_1+2b_1     \\
\lambda_0+\lambda_1-2a_1  & \lambda_0-\lambda_1+2b_1  &  \lambda_0-\lambda_1-2b_1 &  \lambda_0+\lambda_1+2a_1    
\end{bmatrix}.$$ 
Note that $A-JC$ and $A+JC$ are normal matrices with the eigenvalues $z_2, z_3, \overline{z}_2, \overline{z}_3$ and $\lambda_0, \lambda_1, z_1, \overline{z}_1$, respectively. Moreover, since $\lambda_0 + \lambda_1 - 2\sum_{j=1}^3 \vert z_{j} \vert \geq 0$, $A+JC$ is a nonnegative matrix such that each of its entry is not less than the absolute value of the corresponding entry of $A-JC$. Thus $A$ and $C$ are nonnegative matrices. So, the matrix 
$\begin{bmatrix}
A&JCJ\\
C&JAJ
\end{bmatrix}$
is our desired matrix.
\end{proof}
\end{thm}
\medskip
Next, we extend the previous result to the case that there are exactly $(n+1)$ real numbers in our list of $(n+7)$ numbers.
\medskip
\begin{thm} \label{centro10}
Let $\lambda_0 \geq 0 > \lambda_1 \geq \ldots \geq \lambda_n$ be real numbers and $z_1, z_2, z_3$ be complex numbers such that $z_j=a_j+b_{j}i$ where $a_{j} \in \mathbb{R}$ and $b_{j}>0$ for $j=1,2,3$. If $\lambda_0 + \sum_{j=1}^n \lambda_j - 2\sum_{j=1}^3 \vert z_{j} \vert \geq 0$ then there is an $(n+7) \times (n+7)$ normal centrosymmetric nonnegative matrix with the eigenvalues $\lambda_0, \lambda_1, \ldots, \lambda_n, z_1, z_2, z_3, \overline{z}_1, \overline{z}_2, \overline{z}_3.$
\begin{proof}
Suppose that $\lambda_0 + \sum_{j=1}^n \lambda_j - 2\sum_{j=1}^3 \vert z_{j} \vert \geq 0$.
If $n=1$ then the assertion follows from Theorem \ref{centro9}.

Now, suppose that $n=2$ and let $\lambda^{\prime}_0=\lambda_0+\lambda_2$. So, the list of numbers $\lambda^{\prime}_0, \lambda_1, z_1, z_2, z_3$ satisfies the inequality condition of Theorem \ref{centro9}. Then there is an $8 \times 8$ normal centrosymmetric nonnegative matrix $Q$ with the eigenvalues $\lambda^{\prime}_0, \lambda_1, z_1, z_2, z_3, \overline{z}_1, \overline{z}_2, \overline{z}_3$. Set $Q=\begin{bmatrix}
A&JCJ \\
C&JAJ
\end{bmatrix}$ and let $\vtu$ be the unit nonnegative eigenvector of $Q$ corresponding to the Perron eigenvalue $\lambda^{\prime}_0$ such that $J\vtu=\vtu$. So, $\vtu$ can be written as $\vtu=\begin{bmatrix}
\vtx \\ J\vtx
\end{bmatrix}$. Then we extend the matrix $Q$ to the $9 \times 9$ matrix
$\tilde{Q}= \begin{bmatrix}
A            &    JCJ              &   \rho \vtx  \\
C            &    JAJ              &   \rho J\vtx \\
\rho \vtx^T  &    \rho \vtx^TJ     &    0 
\end{bmatrix}$, where $\rho = \sqrt{-\lambda_0 \lambda_2}$. By Theorem \ref{normal1}, $\tilde{Q}$ is a nonnegative normal matrix with the eigenvalues $\lambda_0, \lambda_1, \lambda_2, z_1, z_2, z_3, \overline{z}_1, \overline{z}_2, \overline{z}_3$. By permuting rows and columns, we have that the matrix $\tilde{Q}$ is similar to the matrix 
$\begin{bmatrix}
A            &    \rho \vtx     &   JCJ     \\
\rho \vtx^T  &     0            &   \rho \vtx^TJ \\
C            &    \rho J\vtx    &   JAJ        
\end{bmatrix}$ which is a centrosymmetric matrix and then this matrix is our desired matrix.

Now, we proceed by induction on $n \geq 3$. So assume that the assertion is true for all systems of $\lambda$'s having length less than $n$ and satisfying the condition in the assertion. Define $\lambda^{\prime}_0 = \lambda_0+\lambda_{n}+\lambda_{n-1}$. Then the list of ($n-1$)$+3$ numbers $\lambda^{\prime}_0, \lambda_1, \ldots, \lambda_{n-2}$, $z_1, z_2, z_3$ satisfies the inequality condition in the induction hypothesis. Therefore there is an $(n+5) \times (n+5)$ normal centrosymmetric nonnegative matrix $Q$ with the eigenvalues $\lambda^{\prime}_0, \lambda_1, \ldots, \lambda_{n-2}$, $z_1, z_2, z_3, \overline{z}_1, \overline{z}_2, \overline{z}_3.$ 
Applying Corollary \ref{ourcc3} to the matrices $Q$ and $\begin{bmatrix}
0
\end{bmatrix}$ with $\lambda_{n-1}$ and $\lambda_{n}$, we have that there is an $(n+7) \times (n+7)$ normal centrosymmetric nonnegative matrix with the eigenvalues $\lambda_0, \lambda_1, \ldots, \lambda_n, z_1, z_2, z_3, \overline{z}_1, \overline{z}_2, \overline{z}_3.$
\end{proof}
\end{thm}
\medskip
We now give a result when there are only $4$ real numbers in our list of $(2m+4)$ numbers.
\medskip
\begin{thm} \label{centro13}
Let $\lambda_0 \geq 0 > \lambda_1 \geq \lambda_2 \geq \lambda_3$ be real numbers and 
$z_1, z_2, \ldots, z_{m}$ be complex numbers such that $z_j=a_j + b_ji$ where $a_{j} \in \mathbb{R}$ and $b_j > 0$ for $j=1, \ldots, m$. If $\lambda_0 + \sum_{j=1}^3 \lambda_j - 2 \sum_{j=1}^m \vert z_j \vert \geq 0$, then there is a 
$(2m+4) \times (2m+4)$ normal centrosymmetric nonnegative matrix with all eigenvalues as $\lambda_0$, $\lambda_1$, $\lambda_2$, $\lambda_3$, $z_1$, $\ldots$, $z_m$, $\overline{z}_1$, $\ldots$, $\overline{z}_m$.
\begin{proof}
Suppose that $\lambda_0 + \sum_{j=1}^3 \lambda_j - 2 \sum_{j=1}^m \vert z_j \vert \geq 0$.
If $m=1$ then the assertion follows from Theorem \ref{centro2}. 
If $m \geq 2$ and $m$ is an even number then the assertion follows from Theorem \ref{centro4}.
If $m =3$ then the assertion follows from Theorem \ref{centro10}.
Let $m \geq 5$ and $m=2s+1$ where $s \geq 2$. Define $\lambda^{\prime}_0 = \lambda_0+\lambda_2+\lambda_3-2\sum_{j=1}^3 \vert z_j \vert$ and $\lambda^{\prime \prime}_0 = -\lambda_3 + 2\sum_{j=1}^3 \vert z_j \vert$. Then the list of numbers $\lambda^{\prime}_0, \lambda_1, z_4, z_5, \ldots, z_{m}$ satisfies the inequality condition in Theorem \ref{centro4} and the list of numbers $\lambda^{\prime \prime}_0, \lambda_3, z_1, z_2, z_3$ satisfies the inequality condition in Theorem \ref{centro10}. Therefore there are normal centrosymmetric nonnegative matrices $Q_1$ and $Q_2$ with the eigenvalues $\lambda^{\prime}_0, \lambda_1, z_4, z_5, \ldots, z_{m}, \overline{z}_4, \overline{z}_5, \ldots, \overline{z}_{m}$ and $\lambda^{\prime \prime}_0, \lambda_3, z_1, z_2, z_3, \overline{z}_1, \overline{z}_2, \overline{z}_3$, respectively. Since $Q_1$ and $Q_2$ are of even sizes, we can write $Q_1$ and $Q_2$ as 
$\begin{bmatrix}
A_1 & JC_1J \\
C_1 & JA_1J
\end{bmatrix}$ and
$\begin{bmatrix}
A_2 & JC_2J \\
C_2 & JA_2J
\end{bmatrix}$, respectively. 

If $\lambda^{\prime}_0 \geq \lambda^{\prime \prime}_0$, then we let $\rho = \sqrt{\sigma(\lambda^{\prime}_0-\lambda^{\prime \prime}_0+\sigma)}$ where $\sigma = -\lambda_2-\lambda_3+2\displaystyle \sum_{j=1}^3 \vert z_j \vert$, and if $\lambda^{\prime \prime}_0 > \lambda^{\prime}_0$, then we let $\rho = \sqrt{\sigma(\lambda^{\prime \prime}_0-\lambda^{\prime}_0+\sigma)}$ where $\sigma = \lambda_0 + \lambda_3 - 2\displaystyle \sum_{j=1}^3 \vert z_j \vert$. We combine the matrices $Q_1$ and $Q_2$ to obtain the matrix 
$Q=\begin{bmatrix}
Q_1            & \rho \vtu_0 \vtv^T_0 \\
\rho \vtv_0 \vtu^T_0 & Q_2
\end{bmatrix}$, where $\vtu_0$ is the unit nonnegative eigenvector of $Q_1$ corresponding to the Perron root $\lambda^{\prime}_0$ such that $J\vtu_0=\vtu_0$ and $\vtv_0$ is the unit nonnegative eigenvector of $Q_2$ corresponding to the Perron root $\lambda^{\prime \prime}_0$ such that $J\vtv_0=\vtv_0$. By Theorem \ref{normal1}, the matrix $Q$ is a nonnegative normal matrix with the eigenvalues $\lambda_0, \lambda_1, \ldots, \lambda_n, z_1, \ldots, z_m, \overline{z}_1, \ldots, \overline{z}_m$. Then we can write $\vtu_0=\begin{bmatrix}
\vtx \\ J\vtx
\end{bmatrix}$ and 
$\vtv_0=\begin{bmatrix}
\vty \\ J\vty
\end{bmatrix}$. Finally, we permute rows and columns of the matrix $Q$ to obtain the matrix
$\begin{bmatrix}
A_1               &    \rho \vtx\vty^T    & \rho \vtx\vty^TJ      &  JC_1J \\
\rho \vty \vtx^T  &    A_2                & JC_2J                 &  \rho \vty\vtx^TJ \\
\rho J\vty\vtx^T  &    C_2                & JA_2J                 &  \rho J\vty\vtx^TJ \\
C_1               &    \rho J\vtx\vty^T   & \rho J\vtx\vty^TJ     &  JA_1J
\end{bmatrix}$ which is also a normal centrosymmetric nonnegative matrix with the eigenvalue $\lambda_0$, $\lambda_1$, $\lambda_2$, $\lambda_3$, $z_1$, $\ldots$, $z_m$, $\overline{z}_1$, $\ldots$, $\overline{z}_m$. This is our desired matrix.
\end{proof}
\end{thm}
\medskip

Now, we are ready to give a sufficient condition of the NIEP for normal centrosymmetric matrices where we have a list of $n+2m+1$ complex numbers with $n \geq 3$. This result extends the result of Xu in Theorem \ref{radwan1} when $n \geq 3$.

\medskip
\begin{thm} \label{centromain} 
Let $\lambda_0 \geq 0 > \lambda_1 \geq \lambda_2 \geq \ldots \geq \lambda_n, n \geq 3,$ be real numbers and 
$z_1, z_2, \ldots, z_{m}$ be complex numbers such that $z_j=a_j + b_ji$ where $a_j \in \mathbb{R}$ and $b_j > 0$ for $j=1, \ldots, m$. If $\lambda_0 + \sum_{j=1}^n \lambda_j - 2 \sum_{j=1}^m \vert z_j \vert \geq 0$, then there is an 
$(n+2m+1) \times (n+2m+1)$ normal centrosymmetric nonnegative matrix with the eigenvalues as $\lambda_0, \lambda_1, \ldots, \lambda_n, z_1, \ldots, z_m, \overline{z}_1, \ldots, \overline{z}_m$.
\begin{proof}
Suppose that $\lambda_0 + \sum_{j=1}^n \lambda_j - 2 \sum_{j=1}^m \vert z_j \vert \geq 0$.
If $n=3$ then the assertion follows from Theorem \ref{centro13}.

Now, assume that $n=4$ and let $\lambda^{\prime}_0=\lambda_0+\lambda_4$. Then the list of numbers $\lambda^{\prime}_0$, $\lambda_1$, $\lambda_2$, $\lambda_3$, $z_1$, $z_2$, $\ldots$, $z_m$ satisfies the inequality condition of Theorem \ref{centro13}. Then there is a $(2m+4) \times (2m+4)$ normal centrosymmetric nonnegative matrix $Q$ with the eigenvalues $\lambda^{\prime}_0, \lambda_1, \lambda_2, \lambda_3$, $z_1$, $\ldots$, $z_m$, $\overline{z}_1$, $\ldots$, $\overline{z}_m$. 
Since $Q$ is of even size, we can write $Q$ as $\begin{bmatrix}
A&JCJ \\
C&JAJ
\end{bmatrix}$. Let $\vtu$ be the nonnegative unit eigenvector of $Q$ corresponding to the Perron eigenvalue $\lambda^{\prime}_0$ such that $J\vtu=\vtu$. Then $\vtu$ can be written as $\vtu=\begin{bmatrix}
\vtx \\ J\vtx
\end{bmatrix}$. We extend the matrix $Q$ to the matrix
$\tilde{Q}= \begin{bmatrix}
A         &    JCJ        &   \rho \vtx  \\
C         &    JAJ        &   \rho J\vtx \\
\rho \vtx^T  &    \rho \vtx^TJ  &    0 
\end{bmatrix}$, where $\rho = \sqrt{-\lambda_0 \lambda_4}$. By Theorem \ref{normal1}, $\tilde{Q}$ is a nonnegative normal matrix with the eigenvalues $\lambda_0, \lambda_1, \lambda_2, \lambda_3, \lambda_4, z_1, \ldots, z_m, \overline{z}_1, \ldots, \overline{z}_m$. By permuting rows and columns, the matrix $\tilde{Q}$ is similar to the matrix 
$\begin{bmatrix}
A         &    \rho \vtx     &   JCJ     \\
\rho \vtx^T  &     0         &   \rho \vtx^TJ \\
C         &    \rho J\vtx    &   JAJ        
\end{bmatrix}$ which is a centrosymmetric matrix and it is our desired matrix.

We now proceed by induction on $n \geq 5$. So assume that the assertion is true for all systems of $\lambda$'s having length less than $n$ and satisfying the condition in the assertion. 
Define $\lambda^{\prime}_0 = \lambda_0+\lambda_{n}+\lambda_{n-1}$. Then the list of numbers $\lambda^{\prime}_0, \lambda_1, \ldots, \lambda_{n-2}$, $z_1, \ldots, z_m$ satisfies the condition in the induction hypothesis. Therefore there is a $(2m+n-1) \times (2m+n-1)$ normal centrosymmetric nonnegative matrix $Q$ with the eigenvalues $\lambda^{\prime}_0, \lambda_1, \ldots, \lambda_{n-2}$, $z_1, \ldots, z_m, \overline{z}_1, \ldots, \overline{z}_m$. We extend the matrix $Q$ to the matrix $\tilde{Q}=
\begin{bmatrix}
0            &  \rho \vtu   & -\lambda_{n} \\
\rho \vtu^T     &  Q        & \rho \vtu^TJ    \\
-\lambda_{n} &  \rho J\vtu  & 0
\end{bmatrix}$, where $\vtu$ is the unit eigenvector of the matrix $Q$ corresponding to $\lambda^{\prime}_0$ such that $J\vtu=\vtu$ and $\rho = \sqrt{\dfrac{-(\lambda_0+\lambda_n)(\lambda_{n-1}+\lambda_n)}{2}}$. By Lemma \ref{ourcc2}, we have $\tilde{Q}$ is a $(2m+n+1) \times (2m+n+1)$ normal centrosymmetric nonnegative matrix with the eigenvalues 
$\lambda_0, \lambda_1, \ldots, \lambda_n, z_1, \ldots, z_m, \overline{z}_1, \ldots, \overline{z}_m$. This is our desired matrix.
\end{proof}
\end{thm}
\medskip
Now, we illustrate the construction of a normal centrosymmetric nonnegative matrix using Theorem \ref{centromain}.
\medskip

\begin{exam}
We construct a normal centrosymmetric nonnegative matrix with the eigenvalues $20, -1, -2, -3, 3 \pm 4i, \pm 2i$. Note that this list satisfies the condition in Theorem \ref{centromain}. First, using the proof of Theorem \ref{centro3} with the list of numbers $20+(-2)+(-3)=15, -1, 3+4i, 2i$, we have $$A+JC=
\begin{bmatrix}
7                         &  \dfrac{12+4\sqrt{3}}{3}     & \dfrac{12-4\sqrt{3}}{3}  \\
\dfrac{12-4\sqrt{3}}{3}   &     7                        &  \dfrac{12+4\sqrt{3}}{3} \\
\dfrac{12+4\sqrt{3}}{3}   &   \dfrac{12-4\sqrt{3}}{3}    &   7 
\end{bmatrix}$$ and 
$$A-JC=
\begin{bmatrix}
\dfrac{-1}{3}             &  \dfrac{-1+2\sqrt{3}}{3}     & \dfrac{-1-2\sqrt{3}}{3}  \\
\dfrac{-1-2\sqrt{3}}{3}   &   \dfrac{-1}{3}              &  \dfrac{-1+2\sqrt{3}}{3} \\ 
\dfrac{-1+2\sqrt{3}}{3}   &    \dfrac{-1-2\sqrt{3}}{3}   &    \dfrac{-1}{3}                       
\end{bmatrix}.$$
So, the matrix $\tilde{Q}=\begin{bmatrix}
A & JCJ \\
C & JAJ
\end{bmatrix}$ is a normal centrosymmetric nonnegative matrix with the eigenvalues $15, -1, 3 \pm 4i, \pm 2i$. Finally, we apply Lemma \ref{ourcc2} to $B=\begin{bmatrix}
0
\end{bmatrix}, \rho=\sqrt{\dfrac{85}{2}}$ and $\xi=3$. Then the matrix
$$Q=\begin{bmatrix}
0      & \sqrt{\dfrac{85}{2}}\vtu^T     & 3 \\
\sqrt{\dfrac{85}{2}}\vtu & \tilde{Q}     & \sqrt{\dfrac{85}{2}}\vtu \\
3      & \sqrt{\dfrac{85}{2}} \vtu^T    & 0
\end{bmatrix},$$
where $\vtu^T=\dfrac{1}{\sqrt{6}} \begin{bmatrix}
1&1&1&1&1&1
\end{bmatrix}$ is a unit nonnegative eigenvector of $\tilde{Q}$ corresponding to the Perron root $15$, is a normal centrosymmetric nonnegative matrix with the eigenvalues $20, -1, -2, -3, 3 \pm 4i, \pm 2i$. Explicitly, 
$$Q=\begin{bmatrix}
0&\sqrt{\frac{85}{12}}&\sqrt{\frac{85}{12}}&\sqrt{\frac{85}{12}}&\sqrt{\frac{85}{12}}&\sqrt{\frac{85}{12}}&\sqrt{\frac{85}{12}}&3\\

\sqrt{\frac{85}{12}}& \frac{20}{6} & \frac{11+6\sqrt{3}}{6} & \frac{11-6\sqrt{3}}{6} & \frac{13-2\sqrt{3}}{6} & \frac{13+2\sqrt{3}}{6} & \frac{22}{6} &\sqrt{\frac{85}{12}} \\

\sqrt{\frac{85}{12}}& \frac{11-6\sqrt{3}}{6} & \frac{20}{6} & \frac{11+6\sqrt{3}}{6} & \frac{13+2\sqrt{3}}{6} & \frac{22}{6} & \frac{13-2\sqrt{3}}{6} & \sqrt{\frac{85}{12}} \\

\sqrt{\frac{85}{12}}& \frac{11+6\sqrt{3}}{6} & \frac{11-6\sqrt{3}}{6} & \frac{20}{6} & \frac{22}{6} & \frac{13-2\sqrt{3}}{6} &  \frac{13+2\sqrt{3}}{6} & \sqrt{\frac{85}{12}}  \\

\sqrt{\frac{85}{12}}& \frac{13+2\sqrt{3}}{6} & \frac{13-2\sqrt{3}}{6} & \frac{22}{6} &  \frac{20}{6} & \frac{11-6\sqrt{3}}{6} &   \frac{11+6\sqrt{3}}{6} &     \sqrt{\frac{85}{12}}  \\

\sqrt{\frac{85}{12}} & \frac{13-2\sqrt{3}}{6} & \frac{22}{6} & \frac{13+2\sqrt{3}}{6} &  \frac{11+6\sqrt{3}}{6} & \frac{20}{6} & \frac{11-6\sqrt{3}}{6}  &     \sqrt{\frac{85}{12}}  \\

\sqrt{\frac{85}{12}} & \frac{22}{6} & \frac{13+2\sqrt{3}}{6} & \frac{13-2\sqrt{3}}{6} &  \frac{11-6\sqrt{3}}{6} & \frac{11+6\sqrt{3}}{6} & \frac{20}{6}   &     \sqrt{\frac{85}{12}}  \\

3&\sqrt{\frac{85}{12}}&\sqrt{\frac{85}{12}}&\sqrt{\frac{85}{12}}&\sqrt{\frac{85}{12}}&\sqrt{\frac{85}{12}}&\sqrt{\frac{85}{12}}&0\\
\end{bmatrix}.$$
\end{exam}
\bigskip

\section{Another sufficient condition of the NIEP for normal centrosymmetric matrices}
\bigskip

In this section we present a variance of Theorem \ref{julio5} due to Julio and Soto in the case of normal centrosymmetric nonnegative matrices. To prove this result, we use the following rank-$r$ perturbation theorem due to Julio, Manzaneda and Soto.
\medskip

\begin{thm} \textnormal{(Julio, Manzaneda and Soto \cite{julio2}, 2015)} \label{nrado} 

\noindent Let $A$ be an $n \times n$ normal matrix with the eigenvalues $\lambda_1$, $\lambda_2$, $\ldots$, $\lambda_n$ and for some $r \leq n$, let $\lbrace \vtx_1$, $\vtx_2$, $\ldots$, $\vtx_r \rbrace$ be an orthonormal set of eigenvectors of $A$ corresponding to $\lambda_1, \lambda_2, \ldots, \lambda_r$, respectively. Let $X$ be the $n \times r$ matrix with the $j^{th}$ column $\vtx_j$, let $\Omega = \emph{diag} (\lambda_1, \lambda_2, \ldots, \lambda_r)$, and let $C=\left( c_{ij} \right)$, with $c_{jj} = 0$ for $j=1, \ldots, r$, be any $r \times r$ matrix such that $\Omega+C=B$ is a normal matrix. Then the matrix $A+XCX^*$ is a normal matrix with the eigenvalues $\mu_1, \ldots, \mu_r, \lambda_{r+1}, \ldots, \lambda_{n}$, where $\mu_1, \ldots, \mu_r$ are the eigenvalues of the matrix $B$.
\end{thm}
\medskip
Now, we give the following result.
\medskip
\begin{thm} \label{centromain2} 
Let $\Lambda = \lbrace \lambda_1, \lambda_2, \ldots, \lambda_n \rbrace$ be the set of complex numbers such that $\mathbb{R} \ni \lambda_1 \geq \vert \lambda_j \vert$ for $j=2, \ldots, n$, $\sum_{j=1}^n \lambda_j \geq 0$ and $\overline{\Lambda} = \Lambda$ and let $\omega_1, \ldots, \omega_S$ be nonnegative numbers where $S \leq n$ and $ 0 \leq \omega_k \leq \lambda_1$ for $k=1, \ldots, S.$ Suppose that

\begin{itemize} 
\item[\textnormal{(1)}]  there is a partition $\Lambda_1 \cup \ldots \cup \Lambda_S$ of \ $\lbrace \lambda_1, \ldots, \lambda_n \rbrace$, where for each $j=1,2,\ldots, S$, $\Lambda_j=\lbrace \lambda_{j1}, \lambda_{j2} \ldots \lambda_{jp_j} \rbrace$ with $\lambda_{11}=\lambda_1$, such that at most one set of $\Lambda_j$'s is of odd size and for each $j=1, 2, \ldots, S$, the set $\Gamma_j= \lbrace \omega_j, \lambda_{j2}$, $\ldots$, $\lambda_{jp_j} \rbrace$ is realizable by a normal centrosymmetric nonnegative matrix with the Perron root $\omega_j$, and

\item[\textnormal{(2)}] there is an $S \times S$ normal nonnegative matrix $B$ with all eigenvalues as $\lambda_{11}, \lambda_{21}, \ldots, \lambda_{S1}$ and the diagonal entries as $\omega_S, \omega_{S-1}, \ldots, \omega_1.$
\end{itemize}
Then $\lbrace \lambda_1, \lambda_2, \ldots, \lambda_n \rbrace$ is realizable by a normal centrosymmetric nonnegative matrix.

\begin{proof}

First, we suppose that $\Lambda_j$ is of even size for each $j=1, 2, \ldots, S$. For each $j$, let $\Gamma_j$ be realizable by a normal centrosymmetric nonnegative matrix $Q_j$ of $p_j \times p_j$ size. By Theorem \ref{centroprop1}, we can write 
$Q_j=\begin{bmatrix}
A_j & JC_jJ \\
C_j & JA_jJ
\end{bmatrix}$
, where $A_j$ and $C_j$ are $\frac{p_j}{2} \times \frac{p_j}{2}$ matrices. Then $$\widehat{Q} = 
\begin{bmatrix}
A_S    &                       &      &       &                      &  JC_SJ\\
       & \ddots                &      &       &\reflectbox{$\ddots$} &       \\
       &                       &  A_1 & JC_1J &                      &       \\
       &                       &  C_1 & JA_1J &                      &       \\
       & \reflectbox{$\ddots$} &      &       &\ddots                &       \\
C_S    &                       &      &       &                      &  JA_SJ\\
\end{bmatrix}$$
is a normal centrosymmetric nonnegative matrix with all eigenvalues obtained from all numbers in $\Gamma_j$ for $j=1, \ldots, S$.

By Lemma \ref{centroprop4}, for each $j=1, \ldots, S$ we can find the nonnegative unit eigenvector $\vtx_j= \begin{bmatrix}
\vtv_j \\ J\vtv_j
\end{bmatrix}$ of $Q_j$ corresponding to $\omega_j$.
Then
$$\vtx_1 = \begin{bmatrix} \textbf{0} \\ \vdots \\ \textbf{0} \\ \vtv_1 \\J\vtv_1 \\ \textbf{0} \\ \vdots\\ \textbf{0} \end{bmatrix},
 \vtx_2 = \begin{bmatrix} \textbf{0} \\ \vdots \\ \vtv_2 \\ \textbf{0} \\ \textbf{0} \\J\vtv_2 \\ \vdots\\ \textbf{0} \end{bmatrix},\ldots,
 \vtx_S =\begin{bmatrix} \vtv_S\\ \textbf{0} \\ \vdots \\ \textbf{0} \\ \textbf{0} \\ \vdots \\ \textbf{0} \\ J\vtv_S \end{bmatrix}$$
form an orthonormal set in $\mathbb{R}^n$ and $\vtx_1, \ldots, \vtx_S$ are eigenvectors of $\widehat{Q}$ corresponding to $\omega_1, \ldots, \omega_S$, respectively. Clearly, the $n \times S$ matrix 
$$X=\begin{bmatrix}
\vtx_S & \vtx_{S-1}& \cdots& \vtx_1
\end{bmatrix}$$ 
is a nonnegative matrix with $JX=X$, where $J$ is the $n \times n$ reverse identity matrix, and $\widehat{Q}X=X\Omega$, where $\Omega$ is diag$(\omega_S, \omega_{S-1}, \ldots, \omega_1)$. Therefore $\widehat{Q}+X(B-\Omega)X^*=\widehat{Q}+X(B-\Omega)X^T$ is a normal matrix with the eigenvalues $\lambda_1, \lambda_2, \ldots, \lambda_n$ by Theorem \ref{nrado}. Since the diagonal entries of $B$ are $\omega_S, \omega_{S-1}, \ldots, \omega_1$, $B-\Omega$ is a nonnegative matrix. Thus $\widehat{Q}+X(B-\Omega)X^T$ is also a nonnegative matrix. Moreover, $\widehat{Q}+X(B-\Omega)X^T$ is a centrosymmetric matrix because $\widehat{Q}$ is centrosymmetric and $JX=X$. Therefore $\widehat{Q}+X(B-\Omega)X^T$ is a desired matrix.

If there is a $p \in \lbrace 1, 2, \ldots, S \rbrace $ such that $\vert \Lambda_p \vert$ is odd, then if $\Gamma_p$ is realizable by the nonnegative centrosymmetric matrix $Q_p$, we put $Q_p$ in the center of the matrix $\widehat{Q}$ above and thus we can construct our desired matrix using the same idea as in the above case.
\end{proof}
\end{thm}
\medskip
Our next result is an another variance of Theorem \ref{julio5} in the case of normal centrosymmetric nonnegative matrices.
\medskip
\begin{thm} \label{centromain3} 
Let $\Lambda = \lbrace \lambda_1, \lambda_2, \ldots, \lambda_n \rbrace$ be a set of complex numbers such that $\mathbb{R} \ni \lambda_1 \geq \vert \lambda_j \vert$ for $j=2, \ldots, n$, $\sum_{j=1}^n \lambda_j \geq 0$ and $\overline{\Lambda} = \Lambda$ and let $\omega_1, \ldots, \omega_{m+S}$ be nonnegative numbers where $S$ is a nonnegative integer, $2m+S \leq n$ and $0 \leq \omega_k \leq \lambda_1$ for $k=1, \ldots, m+S.$ Suppose that

\begin{itemize} 
\item[\textnormal{(1)}] there is a partition $(\Lambda_1 \cup \ldots \cup \Lambda_m) \cup (\Lambda_{m+1} \cup \Lambda_{m+2} \cup \ldots \cup \Lambda_{m+S}) \cup (\Lambda_m \cup \ldots \cup \Lambda_1)$ of $\lbrace \lambda_1, \ldots, \lambda_n \rbrace$ where for each $j=1,2,\ldots, m+S$, $\Lambda_j=\lbrace \lambda_{j1}, \lambda_{j2} \ldots \lambda_{jp_j} \rbrace$ such that at most one set of $\Lambda_j$'s is of odd size for $j=m+1, m+2, \ldots, m+S$ and the set $\Gamma_j= \lbrace \omega_j, \lambda_{j2}, \ldots, \lambda_{jp_j} \rbrace$ is realizable by a nonnegative normal matrix for $j \in \lbrace 1,2,\ldots, m \rbrace$ and by a normal centrosymmetric nonnegative matrix for $j \in \lbrace m+1, m+2, \ldots, m+S \rbrace$ with the Perron root $\omega_j$ in every case and

\item[\textnormal{(2)}] there is a $(2m+S) \times (2m+S)$ nonnegative normal matrix
$$B = \begin{bmatrix}
\tilde{A}  & Y^T       & J \tilde{C}J \\
Y          & \tilde{B} & JYJ  \\
\tilde{C}  & JY^TJ     & J \tilde{A} J\\
\end{bmatrix},$$
where $\tilde{B}$ is of size $S\times S$, $B$ has the eigenvalues as $\lambda_{11},$ $\lambda_{21}, \ldots, \lambda_{m1},$ $\lambda_{(m+1)1}$, $\lambda_{(m+2)1}, \ldots, \lambda_{(m+S)1},$ $\lambda_{m1}, \ldots, \lambda_{21},$ $\lambda_{11}$ and the diagonal entries as $\omega_{1},$ $\omega_{2}, \ldots,$ $\omega_{m},$ $\omega_{m+1}$, $\omega_{m+2}, \ldots, \omega_{m+S},$ $\omega_{m}, \ldots, \omega_{2},$ $\omega_{1}$, and $JY=Y$.
\end{itemize}
\noindent Then $\lbrace \lambda_1, \lambda_2, \ldots, \lambda_n \rbrace$ is realizable by a normal centrosymmetric nonnegative matrix.

\medskip
\noindent Remark: If $S=0$, then the condition $(1)$ will be the same except that the partition is of the form $(\Lambda_1 \cup \ldots \cup \Lambda_m) \cup (\Lambda_m \cup \ldots \cup \Lambda_1)$ and we don't need the existence of normal centrosymmetric nonnegative matrices. Also if $S=0$, the matrix $B$ in the condition $(2)$ is of the form
$\begin{bmatrix}
\tilde{A} & J \tilde {C} J \\
\tilde{C} & J \tilde{A} J \\
\end{bmatrix}$ and thus $B$ is a normal centrosymmetric nonnegative matrix in this case.
\begin{proof}
First, we suppose that $S>0$. Now, we consider the case that $\Lambda_j$ is of even size for every $j =m+1, m+2, \ldots, m+S$. For each $j$, let $\Gamma_{j}$ be realizable by a nonnegative normal matrix $P_j$ for $j=1, \ldots, m$ and be realizable by a normal centrosymmetric nonnegative matrix $Q_{j}$ for $j=m+1, \ldots, m+S.$ So, the matrix $Q_j$ is of the form $Q_j = \begin{bmatrix}
A_j & JC_jJ \\
C_j & JA_jJ \\
\end{bmatrix},$ for $j=m+1, \ldots, m+S.$ Let 

$$\widehat{Q}=
\begin{bmatrix}
 P_1&     &        &  &  &  &    \\
    & \ddots &     &  &  &  &    \\
    &     & P_m    &  &  &  &    \\ 
    &     &        & E&  &  &    \\ 
    &     &        &  &JP_{m}J &  &    \\
    &     &        &  &  & \ddots &    \\
    &     &        &  &  &  & JP_{1}J   \\    
\end{bmatrix},$$

\noindent where 

$$E= 
\begin{bmatrix}
A_{m+1} &        &     &         &  & JC_{m+1}J \\
        & \ddots &     &         &\reflectbox{$\ddots$} &  \\
        &        & A_{m+S} & JC_{m+S}J &  &           \\
        &        & C_{m+S} & JA_{m+S}J &  &           \\        
        & \reflectbox{$\ddots$}& &  & \ddots & \\
C_{m+1} &        &     &         &  &  JA_{m+1}J \\
\end{bmatrix}.$$ 

Then $\widehat{Q}$ is a normal centrosymmetric nonnegative matrix with all eigenvalues obtained from $(\Gamma_1 \cup \ldots \cup \Gamma_m) \cup (\Gamma_{m+1} \cup \ldots \cup \Gamma_{m+S}) \cup (\Gamma_m \cup \ldots \cup \Gamma_1)$. Let $\vtu_j$ be nonnegative eigenvectors of $P_j$ corresponding to $\omega_j$ for $j=1, \ldots, m$ and let
$\begin{bmatrix}
\vtv_j \\ 
J\vtv_j \\
\end{bmatrix}$ 
be nonnegative unit eigenvectors of $Q_j$ corresponding to $\omega_j$ for $j=m+1, \ldots, m+S.$ Let 
$X= \begin{bmatrix}
 \tilde{U}&              &             \\
          &  \tilde{X}   &             \\
          &              & J\tilde{U}J    \\
\end{bmatrix}$ where 
$\tilde{U} = 
\begin{bmatrix}
\vtu_1 &        &     \\
    & \ddots &     \\
    &        & \vtu_m \\
\end{bmatrix}$ and

$$\tilde{X}= 
\begin{bmatrix}
\vtv_{m+1} & \cdots  &   0       &    0      \\
    0   & \cdots  & \vdots    & \vdots    \\
    0   & \cdots  &   0       &    0      \\
\vdots  &         & \vtv_{m+S-1} &    0      \\
\vdots  &         &   0       &  \vtv_{m+S}  \\
\vdots  &         &   0       & J\vtv_{m+S}   \\
\vdots  &         & J\vtv_{m+S-1}&    0      \\
    0   & \cdots  &   0       &    0       \\        
    0   & \cdots  & \vdots    &  \vdots     \\
J\vtv_{m+1}& \cdots  &   0       &    0      \\
\end{bmatrix}.$$ 
\noindent Therefore the columns of $X$ consisting of the eigenvectors of $\widehat{Q}$ corresponding to $\omega_1, \ldots, \omega_m,$ $\omega_{m+1},$ $\omega_{m+2}$, $\ldots, \omega_{m+S}$, $\omega_{m}, \ldots, \omega_1$ form an orthonormal set in $\mathbb{R}^n$ and 
$$JX = 
\begin{bmatrix}
            &           &  \tilde{U}J \\
            & \tilde{X} &             \\
J\tilde{U}  &           &             \\
\end{bmatrix} = \begin{bmatrix}
\tilde{U}   &           &              \\
            & \tilde{X} &              \\
            &           &  J\tilde{U}J \\
\end{bmatrix}
\begin{bmatrix}
            &           &    J         \\
            &    I_{S}  &               \\
     J      &           &               \\
\end{bmatrix} = X \begin{bmatrix}
            &           &    J         \\
            &    I_{S}  &               \\
     J      &           &               \\
\end{bmatrix},$$ 
where $I_S$ is an identity matrix of size $S \times S$. Moreover, $\widehat{Q}X = X\Omega$, where 
$\Omega=\begin{bmatrix}
\Omega_1 &          &            \\
         & \Omega_2 &             \\
         &          &  J\Omega_{1}J \\ 
\end{bmatrix},$ $\Omega_1$ = diag($\omega_1, \ldots, \omega_m$), and $\Omega_2$ = diag($\omega_{m+1}, \ldots, \omega_{m+S}$). Therefore $\widehat{Q} + X(B-\Omega)X^*=\widehat{Q} + X(B-\Omega)X^T$ is a nonnegative normal matrix by Theorem \ref{nrado}.
Since the diagonal entries of $B$ are $\omega_1, \ldots, \omega_m$, $\omega_{m+1}$, $\ldots$, $\omega_{m+S}$,  $\omega_{m}$, $\ldots$, $\omega_1$, $B-\Omega$ is a nonnegative matrix and hence $\widehat{Q} + X(B-\Omega)X^T$ is a nonnegative matrix. Since $JY=Y$, $JX(B-\Omega)X^TJ$ which is equal to
$$X \begin{bmatrix}
&&J\\
&I_S&\\
J&&\\
\end{bmatrix}
\begin{bmatrix}
\tilde{A}-\Omega_1  & Y^T       & J \tilde{C}J \\
Y          & \tilde{B}-\Omega_2 & JYJ  \\
\tilde{C}  & JY^TJ     & J (\tilde{A}-\Omega_1) J\\
\end{bmatrix}
\begin{bmatrix}
&&J\\
&I_S&\\
J&&\\
\end{bmatrix}X^T,$$
will be
$X(B-\Omega)X^T$. Therefore 
$$J(\widehat{Q}+X(B-\Omega)X^T)J=J\widehat{Q}J+JX(B-\Omega)X^TJ = \widehat{Q}+X(B-\Omega)X^T.$$ 

\noindent Hence, $\widehat{Q}+X(B-\Omega)X^T$ is a normal centrosymmetric nonnegative matrix with the eigenvalues $\lambda_1, \ldots, \lambda_n.$

Now we consider the case that there is a $p \in \lbrace m+1, m+2, \ldots, m+S \rbrace$ such that $\vert \Lambda_{p} \vert$ is odd. Let $\Gamma_p$ be realizable by a normal centrosymmetric nonnegative matrix $Q_p$ with the Perron root $\omega_p$ corresponding to the nonnegative unit eigenvector $\vtu_p$. Then we put $Q_p$ into the center of the matrix $E$ in the above construction and put the vector $\begin{bmatrix}
0 & \ldots& 0 & \vtu_p^T & 0 &\ldots& 0 \\
\end{bmatrix}^T$ into the $(p-m)^{\text{th}}$ column of $\tilde{X}.$ Therefore we can construct our desired matrix using the same method as in the previous case.

Finaly, if $S=0$, then $X=\begin{bmatrix}
\tilde{U}&             \\
         & J\tilde{U}J \\
\end{bmatrix},$
$JX=XJ$, and $B$ is a normal centrosymmetric nonnegative. So, 
$J(\widehat{Q}+X(B-\Omega )X^T)J$ = $J \widehat{Q}J$+$JXBX^TJ-JX \Omega X^TJ$ = $\widehat{Q}+XJBJX^T+XJ\Omega JX^T$ = $\widehat{Q}+X(B-\Omega)X^T$. Then $\widehat{Q}+X(B-\Omega)X^T$ is our desired matrix.
\end{proof}
\end{thm}
\medskip
\begin{exam} 
We construct a normal centrosymmetric nonnegative matrix with the eigenvalues $20, -1, -1, -2, \pm 2i, \pm 2i, 4 \pm 3i$. First, we note that Theorem \ref{centromain} cannot be applied to this list of numbers because the inequality condition in Theorem \ref{centromain} does not hold. However, this list of numbers satisfies the conditions in Theorem \ref{centromain3} with the partition $\Lambda_1=\lbrace -1, \pm 2i \rbrace$, $\Lambda_2=\lbrace 20, -2, 4 \pm 2i \rbrace$, $\Lambda_1=\lbrace -1, \pm 2i \rbrace$. The set $\Gamma_1=\lbrace 4, \pm 2i \rbrace$ is realizable by a normal nonegative matrix 

\noindent $P_1=\begin{bmatrix}
\frac{4}{3}           &   \frac{4+2\sqrt{3}}{3}    &   \frac{4-2\sqrt{3}}{3}    \\
\frac{4-2\sqrt{3}}{3} &   \frac{4}{3}              &   \frac{4+2\sqrt{3}}{3}  \\   
\frac{4+2\sqrt{3}}{3} &   \frac{4-2\sqrt{3}}{3}    &   \frac{4}{3}              
\end{bmatrix}$
and $\Gamma_2=\lbrace 10, -2, 4 \pm 3i \rbrace$ is realizable by a normal centrosymmetric nonnegative matrix
$Q_2=\begin{bmatrix}
4              & \frac{3}{2}    & \frac{9}{2}   &       0         \\
\frac{9}{2}    &    4           &    0          &    \frac{3}{2}    \\
\frac{3}{2}    &    0           &    4          &    \frac{9}{2} \\
0              &  \frac{9}{2}   &   \frac{3}{2} &       4   
\end{bmatrix}$ (see Theorem \ref{centro1} for the construction of $Q_2$). Now, following the proof of Theorem \ref{centromain3}, we set 

$$\widehat{Q} =\begin{bmatrix}
P_1 &      &       \\
    & Q_2  &       \\
    &      & JP_1J
\end{bmatrix}, \ \Omega=\begin{bmatrix}
4 &    &   \\
  & 10 &   \\
  &    &  4
\end{bmatrix} \text{ and }
X=\begin{bmatrix}
\frac{1}{\sqrt{3}}  &                      &                   \\
\frac{1}{\sqrt{3}}  &                      &                   \\
\frac{1}{\sqrt{3}}  &                      &                   \\
                    &  \frac{1}{2}         &                    \\
                    &  \frac{1}{2}         &                    \\
                    &  \frac{1}{2}         &                    \\
                    &  \frac{1}{2}         &                    \\
                    &                      & \frac{1}{\sqrt{3}} \\
                    &                      & \frac{1}{\sqrt{3}} \\
                    &                      & \frac{1}{\sqrt{3}} \\
                    &                      & \frac{1}{\sqrt{3}} \\
\end{bmatrix}.$$

\noindent Next, we find a normal nonnegative matrix $B$ with the eigenvalues $20$, $-1$, $-1$ and the diagonal entries $4, 10, 4$. One can choose $B$ to be the matrix

\noindent $\begin{bmatrix}
4          &  \sqrt{55}    &  5          \\
\sqrt{55}  &  10           &  \sqrt{55}  \\
5          &  \sqrt{55}    &  4 
\end{bmatrix}.$ Therefore $Q = \widehat{Q}+X(B-\Omega) X^T$ is a normal centrosymmetric nonnegative matrix with the eigenvalues $20, -1, -1, -2, \pm 2i, \pm 2i, 4 \pm 3i$. Explicitly 

$$Q=\begin{bmatrix}
\frac{4}{3}           &  \frac{4+2\sqrt{3}}{3}  & \frac{4-2\sqrt{3}}{3}  & \frac{\sqrt{55}}{2\sqrt{3}} & \frac{\sqrt{55}}{2\sqrt{3}} & \frac{\sqrt{55}}{2\sqrt{3}} & \frac{\sqrt{55}}{2\sqrt{3}} & \frac{5}{3} & \frac{5}{3} & \frac{5}{3}  \\

\frac{4-2\sqrt{3}}{3} &   \frac{4}{3}              &   \frac{4+2\sqrt{3}}{3} & \frac{\sqrt{55}}{2\sqrt{3}} & \frac{\sqrt{55}}{2\sqrt{3}} & \frac{\sqrt{55}}{2\sqrt{3}} & \frac{\sqrt{55}}{2\sqrt{3}} & \frac{5}{3} & \frac{5}{3} & \frac{5}{3}  \\

\frac{4+2\sqrt{3}}{3} &   \frac{4-2\sqrt{3}}{3}    &   \frac{4}{3}  & \frac{\sqrt{55}}{2\sqrt{3}} & \frac{\sqrt{55}}{2\sqrt{3}} & \frac{\sqrt{55}}{2\sqrt{3}} & \frac{\sqrt{55}}{2\sqrt{3}} & \frac{5}{3} & \frac{5}{3} & \frac{5}{3} \\

\frac{\sqrt{55}}{2\sqrt{3}}& \frac{\sqrt{55}}{2\sqrt{3}} & \frac{\sqrt{55}}{2\sqrt{3}} &
4 & \frac{3}{2} & \frac{9}{2} & 0 &
\frac{\sqrt{55}}{2\sqrt{3}}& \frac{\sqrt{55}}{2\sqrt{3}} & \frac{\sqrt{55}}{2\sqrt{3}} \\

\frac{\sqrt{55}}{2\sqrt{3}}& \frac{\sqrt{55}}{2\sqrt{3}} & \frac{\sqrt{55}}{2\sqrt{3}} &
\frac{9}{2} & 4 & 0 & \frac{3}{2} &
\frac{\sqrt{55}}{2\sqrt{3}}& \frac{\sqrt{55}}{2\sqrt{3}} & \frac{\sqrt{55}}{2\sqrt{3}} \\

\frac{\sqrt{55}}{2\sqrt{3}}& \frac{\sqrt{55}}{2\sqrt{3}} & \frac{\sqrt{55}}{2\sqrt{3}} &
\frac{3}{2} & 0 & 4 & \frac{9}{2} &
\frac{\sqrt{55}}{2\sqrt{3}}& \frac{\sqrt{55}}{2\sqrt{3}} & \frac{\sqrt{55}}{2\sqrt{3}} \\

\frac{\sqrt{55}}{2\sqrt{3}}& \frac{\sqrt{55}}{2\sqrt{3}} & \frac{\sqrt{55}}{2\sqrt{3}} &
0 & \frac{9}{2} & \frac{3}{2} & 4 &
\frac{\sqrt{55}}{2\sqrt{3}}& \frac{\sqrt{55}}{2\sqrt{3}} & \frac{\sqrt{55}}{2\sqrt{3}} \\

\frac{5}{3} & \frac{5}{3} & \frac{5}{3} & 
\frac{\sqrt{55}}{2\sqrt{3}}& \frac{\sqrt{55}}{2\sqrt{3}} & \frac{\sqrt{55}}{2\sqrt{3}} & \frac{\sqrt{55}}{2\sqrt{3}} &
\frac{4}{3}           &  \frac{4-2\sqrt{3}}{3}  & \frac{4+2\sqrt{3}}{3} \\

\frac{5}{3} & \frac{5}{3} & \frac{5}{3} & 
\frac{\sqrt{55}}{2\sqrt{3}}& \frac{\sqrt{55}}{2\sqrt{3}} & \frac{\sqrt{55}}{2\sqrt{3}} & \frac{\sqrt{55}}{2\sqrt{3}} &
\frac{4+2\sqrt{3}}{3}  & \frac{4}{3} & \frac{4-2\sqrt{3}}{3} \\

\frac{5}{3} & \frac{5}{3} & \frac{5}{3} & 
\frac{\sqrt{55}}{2\sqrt{3}}& \frac{\sqrt{55}}{2\sqrt{3}} & \frac{\sqrt{55}}{2\sqrt{3}} & \frac{\sqrt{55}}{2\sqrt{3}} &
\frac{4-2\sqrt{3}}{3}  &  \frac{4+2\sqrt{3}}{3}  & \frac{4}{3} \\
\end{bmatrix}.$$

\end{exam}
\bigskip





\end{document}